\documentclass{amsart}
\usepackage{amssymb}
\usepackage{latexsym}

\newcommand{\R}{{\mathbb R}}

\newcommand{\N}{{\mathbb N}}

\newcommand{\Z}{{\mathbb Z}}

\newcommand{\Sub}{\mbox{Sub}}
\newcommand{\CalW}{{\mathcal{W}(S)}}
\newcommand{\OmegaS}{\Omega(S)}
\newcommand{\Oomega}{(\Omega(S),T)}

\newtheorem{theorem}{Theorem}
\newtheorem{lemma}{Lemma}[section]
\newtheorem{prop}[lemma]{Proposition}
\newtheorem{coro}{Corollary}
\newtheorem{definition}[lemma]{Definition}
\newtheorem{remark}{Remark}

\begin{document}
\title[Substitutions Over Finite Alphabets]{Substitution Dynamical Systems:
Characterization of Linear Repetitivity\\and Applications}
\author[D.~Damanik, D.~Lenz]{David Damanik$\,^{1}$ and Daniel
  Lenz$\,^{2}$}
\thanks{D.\ D.\ was supported in part by NSF Grant No.~DMS--0227289\\
\indent D.\ L.\ was supported in part by the DFG within the priority program Quasikristalle.}

\maketitle
\vspace{0.3cm}
\noindent
$^1$ Department of Mathematics 253--37, California Institute of Technology,
Pasadena, CA 91125, U.S.A., E-Mail: damanik@its.caltech.edu\\[0.2cm]
$^2$ Fakult\"at f\"ur Mathematik, TU Chemnitz, D-09107 Chemnitz, Germany, E-Mail: dlenz@mathematik.tu-chemnitz.de\\[0.3cm]

\begin{abstract}   We consider  dynamical systems arising from substitutions over a finite alphabet. We prove that such a system is linearly repetitive if and only if it is minimal.   Based on this characterization we extend various results from primitive substitutions to minimal substitutions. This includes applications to random Schr\"odinger operators and to number theory.
\end{abstract}

\section{Introduction}  \label{Introduction}
This paper deals with a special class of low complexity subshifts
over finite alphabets, viz.\ subshifts associated to
substitutions.

Subshifts over finite alphabets  play a role in various branches
of mathematics, physics, and computer science. Low complexity or
intermediate disorder has been a particular focus of research in
recent years. This has been even more the case due to the
discovery by Shechtman et al.\ of special solids \cite{SBGC},
later called quasicrystals, which exhibit this form of disorder
\cite{BM, Jan, Sen}. Subshifts associated to substitutions and in
particular to primitive substitutions are foremost among the
models of low complexity subshifts \cite{Lot1, Lot2, Que}.

With the recent work of Durand \cite{Du2} and Lagarias and
Pleasants \cite{LP} it became apparent that a key feature to be
studied in low complexity subshifts (and their higher-dimensional
analogue) is linear repetitivity or linear recurrence. It is known
that subshifts associated to primitive substitutions are linearly
repetitive \cite{Sol}. Thus, it is natural to ask:

\begin{itemize}
\item[(Q)]  Which substitution dynamical systems are linearly repetitive?
\end{itemize}
The main result of  the paper answers this question. Namely, we
show that a substitution dynamical system is linearly repetitive
if and only if it is minimal, which in turn is the case if and
only if one letter (not belonging to a particular subset of the
alphabet) appears with bounded gaps. This not only characterizes
linear repetitivity but also gives an easy to handle  condition to
verify this feature. This characterization and its proof give very
direct methods to
\begin{itemize}
\item[(E)] extend results from the framework of primitive
substitutions to the framework of minimal substitutions.
\end{itemize}

We illustrate this extension process with two types of examples.
The first type is concerned with the spectral theory of certain
Schr\"odinger operators. The second example deals with number
theory. Details will be discussed in the corresponding sections.

The paper is organized as follows: In Section \ref{Notation} we
introduce the necessary notation and state our main result,
answering (Q). This result is then proved in
Section~\ref{Linearly}. The following two sections give examples
for (E). Section \ref{Theassociated} is devoted to a study of
Schr\"odinger operators associated to minimal substitutions. An
application to number theory is discussed in Section \ref{Fixed}.
Finally, in Section \ref{Unique} we study the unique decomposition
property for a special class of nonprimitive substitutions.

\section{Notation and Statement of the Main Result}\label{Notation}
In this section we introduce the necessary notation and present
our main result.

Let $A$ be a finite subset of $\R$, called the alphabet. The
elements of $A$ will be called letters. In the sequel we will use
freely notions from combinatorics on words (see, e.g.,
\cite{Lot1,Lot2}). In particular, the elements of the free monoid
$A^\ast$ over $A$ will be called words. The length of a word is
the number of its letters; the number of occurrences of $v\in
A^\ast$ in $w\in A^\ast$ will be denoted by $\#_v (w)$. Moreover,
for a word $u$ over $A$, we let $\Sub(u)$ denote the set of
subwords of $u$.

We can equip $A$ with discrete topology and $A^{\Z}$ with product
topology. A pair $(\Omega,T)$ is then called a subshift over $A$
if $\Omega$ is a closed subset of $A^{\Z}$ which is invariant
under $T : A^{\Z} \longrightarrow A^{\Z}$, $(T u) (n)\equiv
u(n+1)$. To a subshift $(\Omega,T)$ belongs the  set $\mathcal{W}
(\Omega)$ of finite words given by $ \mathcal{W} (\Omega) \equiv
\cup_{\omega\in \Omega} \Sub(\omega)$. A word $v \in \mathcal{W}
(\Omega) $ is said to occur with bounded gaps if there exists an
$L_v >0$ such that every $w\in \mathcal{W} (\Omega)$ with $|w|\geq
L_v$ contains a copy of $v$. By standard arguments $(\Omega,T)$ is
minimal (i.e., each orbit is dense) if and only if every $v\in
\mathcal{W} (\Omega)$ occurs with bounded gaps.  A particular
strengthening of minimality is thus the condition of linear
repetitivity given as follows: The  system $(\Omega,T)$ is said to
be linearly repetitive if there exists a constant $C_{{\rm LR}}$
with

\begin{equation}\label{lr}
\#_v (w)\geq 1 \;\:\mbox{whenever} \:\; |w|\geq C_{{\rm LR}} |v|
\end{equation}
for $v,w\in \mathcal{W}(\Omega)$.

A special way to generate subshifts is given as follows. Consider
a map $S : A\longrightarrow A^\ast$. By definition of $A^\ast$,
$S$ can be uniquely extended to a morphism
$$
S: A^\ast \longrightarrow A^\ast \;\:\mbox{by setting}\;\:
S(a_1\ldots a_n)\equiv S(a_1)\ldots S(a_n),
$$
for arbitrary $a_j\in A$, $j=1\ldots n$. To such an $S$, we can
associate the set $\CalW\subset A^\ast$ given by
$$
\CalW \equiv \{ w\in A^\ast : w\in \Sub(S^n (a))\;\:\mbox{for
suitable $a\in A$ and $n\in \N_0$ } \},
$$
and the (possibly empty) set $\OmegaS\subset A^{\Z}$ given by
$$\OmegaS \equiv \{\omega\in A^{\Z} : \Sub(\omega)\subset \CalW\}.$$
For $\OmegaS$ to be nonempty, it is necessary and sufficient that
there exists an $e\in A$ with
\begin{equation}\label{length}
|S^n (e)|\longrightarrow \infty, n\longrightarrow \infty.
\end{equation}
Without loss of generality we can then assume (after possibly
removing some letters from $A$) that
\begin{equation}\label{occurrence}
\mbox{for all $a\in A$ there exists an $n\in \N$ with $ \#_a (S^n (e))\geq 1$.}
\end{equation}
Finally, one needs that finite and infinite words associated to $S$ are compatible in the sense that
\begin{equation}\label{compatible}
\CalW = \mathcal{W} (\Omega(S)). 
\end{equation}

\begin{definition} $\Oomega$ is called a substitution dynamical system if \eqref{length}, \eqref{occurrence} and \eqref{compatible} hold.
\end{definition}

\begin{remark}{\rm 
(a) The conditions \eqref{length}, \eqref{occurrence} and \eqref{compatible} are clearly met if the powers $S^n (e)$ converge to a fixpoint of $S$ for $n\longrightarrow \infty$ (i.e. $|S^n (e)|\to \infty$, $n\to \infty$, and $S^n (e)$ is a prefix of $S^{n+1} (e) $ for every $n\in \N$). This is the usual  way to generate a substitution dynamical system. \\[1mm]
(b) As brought to our attention by H. Yuasa,  conditions \eqref{length} and \eqref{occurrence} together do not imply \eqref{compatible}, as can be seen by considering $S : \{0,1\} \longrightarrow \{0,1\}$,  $0\mapsto 1 0$, $1\mapsto 1$. This is an example of what is called ``quasi-primitive substitution'' in \cite{Yua}, where a  further study of such  substitutions can be found.}
\end{remark}

As mentioned in the introduction, a special class of substitution dynamical  systems
$\Oomega$ known to be linearly repetitive are those coming from
primitive $S$. Here, $S$ is called primitive if there exists an
$r\in \N$ with $\#_a (S^r(b))\geq 1$ for arbitrary $a,b\in A$. Our
main result characterizes all $S$ with linearly repetitive
$\Oomega$.

\begin{theorem}\label{main}
Let $\Oomega$ be a substitution dynamical system. Then the
following are equivalent:
\begin{itemize}
\item[{\rm (i)}] There exists an $e\in A$ satisfying
\eqref{length} and  \eqref{occurrence} which occurs with bounded
gaps. \item[{\rm (ii)}] $\Oomega$ is minimal. \item[{\rm (iii)}]
$\Oomega$ is linearly repetitive.
\end{itemize}
\end{theorem}

\begin{remark}{\rm
(a) Condition (i) can be easily checked in many concrete cases.
For example, it can immediately be seen to be satisfied in the
examples of de~Oliveira and Lima \cite{OL}. Thus, their examples
are linearly repetitive and the whole theory developed below
applies.\\[1mm]
(b) As primitive substitutions can easily be seen to satisfy (i), the theorem contains
the well known result (see, e.g., \cite{Sol}) that these systems are linearly
repetitive.\\[1mm]
(c) Existence of an arbitrary $a\in A$ occurring with bounded gaps is not sufficient
for minimality, as can be seen by considering the example $S : \{0,1\} \longrightarrow
\{0,1\}^\ast$, $0\mapsto 1 0 1$ and $1\mapsto 1$.
}
\end{remark}

As discussed above, linear repetitivity implies minimality.
Moreover, it also implies unique ergodicity as shown by Durand
\cite{Du2} (see \cite{Len1} for a different proof as well). Thus,
we obtain the following corollary.

\begin{coro} \label{ue}
Let $\Oomega$ be as above. If $e\in A$ satisfying \eqref{length}
and \eqref{occurrence} occurs with bounded gaps, then $\Oomega$ is
uniquely ergodic and minimal.
\end{coro}

\begin{remark}{\rm  (a) A different way of stating the corollary would be to
say that for a substitution dynamical system, minimality is
equivalent to strict ergodicity. Note that unique ergodicity is
not sufficient for minimality, as can be
seen by considering the example in (c) of Remark~1.\\[1mm]
(b) By (a) of Remark~1, the corollary applies to the examples of
\cite{OL} and we recover their Proposition~1. }
\end{remark}

\section{Linearly Repetitive Substitutions} \label{Linearly}
This section is devoted to a proof of Theorem \ref{main}. It turns
out that the  hard part in the proof is the implication ${\rm
(ii)} \Longrightarrow {\rm (iii)}$. Its proof will be split into
several parts. The key issue is to study the growth of $|S^n(a)|$
for $n\to \infty$ and $a\in A$. This will be done by relating $S :
A \longrightarrow A^\ast$ to a suitable other substitution
$\widetilde{S} : C\longrightarrow C^\ast$, which can be shown to
be primitive if (ii) is satisfied.  As growth properties are well
known for primitive substitutions, we obtain the desired results
by comparing the growth behavior of $S$ and $\widetilde{S}$.

For notational convenience, we will say that $\Oomega$ satisfies
the bounded gap condition (for $e\in A$) if
\begin{itemize}
\item[(BG)] the letter $e\in A$ satisfies \eqref{length} and
\eqref{occurrence} and occurs with bounded gaps.
\end{itemize}

Our first result gives some immediate consequence of (BG).

\begin{lemma}\label{minimal} Let $\Oomega$ satisfy {\rm (BG)}. Then $\Oomega$ is minimal.
\end{lemma}
\begin{proof} It suffices to show that  $S^n (e)$ occurs with bounded gaps for arbitrary
but fixed $n\in \N$.  Set $M\equiv \max\{ |S^n (a)| : a\in A\}$.
By (BG), there exists $\kappa>0$ such that every word in $\CalW$
with length exceeding $\kappa$ contains $e$. Consider an arbitrary
$w\in \CalW$ with $|w|\geq (3 + \kappa) M$. Now, $w$  is contained
in $S^n (a_1\ldots a_s)$ with suitable $a_j\in A$, $j=1,\ldots ,
s$ with $a_1\ldots a_s\in \CalW$. By assumption on $|w|$ and
definition of $\kappa$, we infer that  $w$ contains  $S^n (a_l)
S^n (a_{l+1}) \ldots S^n (a_{l+\kappa})$ for a suitable $l$. By
definition of $\kappa$, we infer that $w$ contains $S^n (e)$ and
the proof is finished.
\end{proof}

We will now introduce the substitution $\widetilde{S}$.  Set
$$B\equiv \{ a\in A :  \limsup_{n\to\infty} |S^n (a)|<\infty\}$$
and
$$ C\equiv A \setminus B.$$
Note that $S$ maps $B^\ast$ into itself. We define
$$\widetilde{S} : C \longrightarrow C^\ast, \:\;\mbox{by}\;\: \widetilde{S}(x)
\equiv \widetilde{S(x) },$$
where for an arbitrary word $w\in
A^\ast$ we define $\widetilde{w}$ to be the word obtained from $w$
by removing every element of $B$.   As $B$ is invariant under $S$,
we infer that
$$ \widetilde{S}^n (\widetilde{x})= \widetilde{ S^n(x)}$$
for arbitrary $x\in \CalW$ and $n\in \N$.  This will be used
repeatedly in the sequel. Our next aim is to show that
$\widetilde{S}$ is primitive if (BG) is satisfied. We need two
preparatory lemmas.

\begin{lemma}\label{eins} Let $\Oomega$ satisfy {\rm (BG)}. Then,
the following are equivalent for $a \in A$:
\begin{itemize}
\item[{\rm (i)}] $|S^n (a)|\longrightarrow \infty$, $n\to \infty$.
\item[{\rm (ii)}] $e$ is contained in $S^k (a)$ for a suitable $k\in \N$.
\end{itemize}
\end{lemma}

\begin{proof} (i) $\Longrightarrow$ (ii). This is clear as $e$ occurs with bounded gaps by (BG). \\
(ii)$\Longrightarrow$ (i). By (ii),  $S^n (e) $ is a subword of $S^{k+n} (a)$ for every $n\in \N$. Now, (i) follows as $e$ satisfies \eqref{length} by (BG).
\end{proof}

\begin{lemma}\label{drei} Let $\Oomega$ satisfy {\rm (BG)}. Then, there exists $m\in \N$ such that $S^n (e)$ contains every letter of $A$ for every $n\geq m$.
\end{lemma}
\begin{proof} By (BG) and \eqref{length}, there exists $r\in \N$  such that $S^n (e) $ contains $e$  whenever $n\geq r$. By \eqref{occurrence}, for every $a\in A$, there exists $n(a)\in \N$ such that $S^{n(a)}(e)$ contains $a$. Then $m= r + \sum_{a\in A} n(a)$ has the desired properties.
\end{proof}

We can now show that $\widetilde{S}$ is primitive, if (BG) holds.

\begin{lemma}\label{primitive} Let $\Oomega$ satisfy {\rm (BG)}. Then $\widetilde{S} : C\longrightarrow C^\ast$ is primitive.
\end{lemma}
\begin{proof} For $c\in C$, we can choose by Lemma \ref{eins} a number $n(c)\in \N$ such that $S^{n(c)} (c)$ contains $e$. Moreover, by Lemma~\ref{drei}, there exists $m$ such that  $S^n (e)$ contains every letter of $A$ whenever $n \geq m$. Let $N\equiv m + \sum_{c\in C} n(c)$. Then, for every $c\in C$, $S^N (c)$ contains every letter of $A$. In particular, for each $c \in C$, the word $ \widetilde{S}^N (c) =\widetilde{ S^n (c)}$ contains every letter of $C$ and primitivity of $ \widetilde{S}$ is proved.
\end{proof}

As $\widetilde{S}$ is primitive, for $c \in C$, the behavior of
$|\widetilde{S}^n (c)|$ for large $n\in \N$ is rather explicit.
The next lemma allows us to compare this behavior with the
behavior of  $|S^n (c)|$.

\begin{lemma}\label{compare}
Let $\Oomega$ satisfy {\rm (BG)}. There exist constants $L>0$ and $N\in \N$ with
$$
\frac{1}{L} \leq \frac{ | \widetilde{S}^n (\widetilde{v})| }{ |S^n (v)|} \leq 1
$$
for arbitrary $n\geq N$ and $v\in \CalW$ containing at least one
letter of $C$.
\end{lemma}

\begin{proof} The inequality $|\widetilde{S}^n (\widetilde{v})| \leq |S^n (v)|$ is obvious. To show the other inequality, note that, by (BG), there exists a constant $\kappa$ such that every word with length exceeding $\kappa$ contains a copy of $e$. This implies, $|\widetilde{v} |\geq \kappa^{-1}  |v| -2$ for arbitrary $v\in \CalW$. Applying this inequality to $\widetilde{S}^n (\widetilde{v}) = \widetilde{S^n (v)}$, we find
$$
|\widetilde{S}^n (\widetilde{v})| \geq \frac{1}{\kappa} |S^n (v)| - 2.
$$
By definition of $C$, there exists $N\in\N$ with
$$
|S^n(c)|\geq 4 \kappa \;\:\mbox{for all $c\in C$ and $n\geq N$.}
$$
Thus,
$$ |\widetilde{S}^n (\widetilde{v})|= |\widetilde{S^n (v)}|
\geq \frac{1}{\kappa} |S^n (v)| - 2 \geq \frac{1}{2\kappa} |S^n
(v)|
$$
for arbitrary $n\geq N$ and $v\in \CalW$  containing at
least one letter of $C$.\end{proof}

The key technical result in this section is the following proposition.

\begin{prop}\label{growth} Let $\Oomega$ satisfy {\rm (BG)}. Let $V$ be a finite subset of $\CalW$ all of whose elements contain at least one letter of $C$. Then there exist $\theta>0$ and $\lambda(V),\rho(V)  >0$ with
$$ \lambda(V) \theta^n \leq |S^n (v) |\leq \rho(V) \theta^n$$
for arbitrary $n\in \N$ and $v\in V$.
\end{prop}
\begin{proof} Set $\widetilde{V}\equiv \{\widetilde{v} : v\in V\}$. As $\widetilde{S}$ is primitive by Lemma \ref{primitive}, there exist $\theta>0$ and constants $\kappa_1,\kappa_2>0$ with $\kappa_1 \leq \theta^{-n}| \widetilde{S}^n (c)| \leq \kappa_2$ for arbitrary  $c\in C$ and  $n\in \N$. As $V$ is finite and every $v\in V$ contains at least one letter of $C$, this shows existence of  $\nu_1,\nu_2>0$  with
$$
\nu_1 \leq \frac{| \widetilde{S}^n (\widetilde{v})|  }{ \theta^n  } \leq \nu_2
$$
for every $\widetilde{v} \in \widetilde{V}$. Therefore, by Lemma \ref{compare}, there exist $\mu_1,\mu_2>0$ and $N\in \N$ with
$$
\mu_1 \theta^n \leq |S^n (v) |\leq \mu_2 \theta^n
$$
for arbitrary $n\geq N $ and $v\in V$.  Adjusting the constants to
fit in the remaining finitely many cases, we conclude the proof.
\end{proof}

With these preparations out of the way, we are now ready to prove Theorem~\ref{main}.

\begin{proof}[Proof of Theorem~\ref{main}.] The implications
(iii) $\Longrightarrow  $ (ii)$ \Longrightarrow $ (i) are obvious.
The implication (i) $\Longrightarrow $ (ii) is given in Lemma
\ref{minimal}.

It remains to prove (ii)$\Longrightarrow $ (iii). We will use the
notion of return word introduced recently by Durand \cite{Du1}.
Recall that $x\in \CalW$ is called a return word of $v\in \CalW$
if $xv\in \CalW$, $xv$ begins with $v$, and $\#_v (xv) =2 $. Let
$e\in A$ satisfying (BG) be fixed. Such an $e$ exists by
minimality. Let $V$ be the set of return words of $e$. As $e$
satisfies (BG), $V$ is a finite set.  Let $U\equiv \{ z_1 z_2 :
z_1, z_2\in V, z_1 z_2\in \CalW\}$. As $V$ is finite, so is $U$.
By the minimality assumption (ii), there exists $G > 0$ such that
every word in $\CalW$ with length exceeding $G$ contains every
word of $U$. By Proposition~\ref{growth}, there exist $\theta,
\lambda(V), \rho(V)>0$ with
\begin{equation}\label{gro}
\lambda(V) \theta^n \leq |S^n (v) |\leq \rho(V) \theta^n
\end{equation}
for all $v\in V$ and $n\in \N$. Define
$$C_{{\rm LR}} \equiv  (3 + G) \theta \rho(V) \lambda(V)^{-1}. $$
We will show linear repetitivity of $\Oomega$ with this constant.
Thus, let $w\in \CalW$ be given and consider an arbitrary   $u\in
\CalW$ with  $|u|\geq C_{{\rm LR}} |w|$. We have to show that $u$
contains a copy of $w$. To do so, we will show the following:

\begin{itemize}
\item $w$ is contained in $S^n( z_0 )$ with suitable $n\in \N$ and  $z_0\in U$,
\item $u$ contains all words of the form $S^n (z)$ with $z\in U$.
\end{itemize}
Here are the details: Let $n\in \N$ be given  with
\begin{equation}
\lambda(V) \theta^{n-1} \leq |w| < \lambda(V) \theta^n.
\end{equation}
Combining this inequality with \eqref{gro}, we see that
\begin{equation}\label{stern}
|w| \leq |S^n (v)| \;\:\mbox{for every $v\in V$}.
\end{equation}
Apparently, we can choose $x= e y e\in \CalW$ such that  $w$ is a
subword of $S^n (e y e)$. Partitioning $e y e$ according to
occurrences of $e$, we can write $e y e = x_1\ldots x_k e$ with
$x_j\in V$, $j=1,\ldots ,k$. By \eqref{stern}, and since $w$ is a
subword of $S^n (x_1)\ldots S^n (x_k) S^n (e)$, we then infer that
$w$ is in fact a subword of $S^n (z_1 z_2)$ with $z_1,z_2\in V$
and $z_1 z_2\in U$. Let us now turn our attention to $u$. By
$|u|\geq C_{{\rm LR}} |w|$, $|w|\geq  \lambda(V) \theta^{n-1}$ and
the definition of $C_{{\rm LR}}$,  we infer
\begin{equation}\label{sternstern}
 (3+G) \rho(V) \theta^{n} \leq |u|.
\end{equation}
Of course, as discussed above for $w$, we can also exhibit $u$ as
a subword of $S^n (x_1 \ldots x_k ) S^n(e)$ with $x_j\in V$,
$j=1,\ldots, k$ and $x_1\ldots x_k\in \CalW$. By \eqref{gro} and
\eqref{sternstern}, we then  conclude that $u$ must contain  a
word of the form $S^n (v)$ with $|v|\geq G$. By definition of $G$,
the word $v$ then contains $z_1 z_2$. Thus, $u$ contains $S^n (z_1
z_2) $ which contains $w$. This finishes the proof.
\end{proof}

\section{The Associated Schr\"odinger Operators} \label{Theassociated}
In this section we discuss applications to Schr\"odinger operators.

Recall that to a given subshift $(\Omega,T)$ over $A\subset \R$,
we can associate the family  $(H_\omega)_{\omega\in \Omega}$ of
selfadjoint operators $H_\omega : \ell^2 (\Z) \longrightarrow
\ell^2 (\Z)$, $\omega\in \Omega$ acting  by
\begin{equation}\label{family}
(H_\omega u) (n)\equiv u(n+1) + u(n-1) + \omega (n) u(n).
\end{equation}
Assume furthermore that $(\Omega,T)$ is minimal, uniquely ergodic,
and aperiodic (i.e., $T^k\omega \neq \omega $ for every $k\neq 0$
and $\omega \in \Omega$). Denote the unique $T$-invariant
probability measure by $\mu$. Such operators have attracted a lot
of attention in recent years  (see, e.g., \cite{Dam,Sut3} for
reviews and below for literature concerning special classes). They
arise in the quantum mechanical treatment of (one-dimensional)
quasicrystals. The theoretical study of physical features (e.g.,
conductance) is accordingly performed by investigating the
spectral theory of such families. It turns out that the spectral
theory of these families is rather interesting. Namely, they
exhibit features such as
\begin{itemize}
\item purely singular continuous spectrum,
\item Cantor spectrum of Lebesgue measure zero,
\item anomalous transport.
\end{itemize}
In the study of these and related properties, two classes of
examples have received particular attention. These are Sturmian
models (and more generally circle maps)
\cite{BIST,Dam3,DL,DL1,DL2,DL4,Kam,Sut} and operators associated
to primitive substitutions \cite{Bel,BBG,BG,Dam4,Dam5, Len2,LTWW}.
The aim of this section is to extend the theory from primitive
substitutions to minimal substitutions, thereby giving a precise
sense to (E) in this case.

For linearly repetitive systems it was recently shown by one of
the authors \cite{Len2} that their spectrum is a Cantor set if
they are not periodic. Thus, we obtain the following  result as an
immediate  corollary to Theorem~\ref{main} above and Corollary~2.2
of \cite{Len2}.

\begin{theorem}\label{zeromeasure} Let $\Oomega$ be an aperiodic
minimal substitution dynamical system. Then, there exists a Cantor
set $\Sigma\subset \R$ of Lebesgue measure zero with
$\sigma(H_\omega)=\Sigma$ for every $\omega\in \Omega$, where
$\sigma(H_\omega)$ denotes the spectrum of $H_\omega$.
\end{theorem}

\begin{remark}{\rm (a) If the subshift is periodic, it is well known that the spectrum
is a finite union of (nondegenerate) closed intervals. Hence, in
this case it is neither a Cantor set nor does it have Lebesgue measure zero.\\[1mm]
(b)  This theorem contains, in particular, the corresponding result for primitive
substitutions obtained in \cite{Len2} (see also \cite{LTWW} for a different proof).\\[1mm]
(c) The theorem covers all the examples discussed in \cite{OL}.
}
\end{remark}

Next we state our result on singular continuous spectrum.

\begin{theorem}\label{singularcont} Let $\Oomega$ be a minimal substitution dynamical
system.  If there exists $u\in \CalW$ starting with $e\in C$ such that $uuue\in\CalW$,
then the operators $(H_\omega)$ have purely singular continuous spectrum for $\mu$-almost
every $\omega\in \Omega$.
\end{theorem}

\begin{remark}{\rm This also covers all the examples studied by de~Oliveira and Lima in
\cite{OL}.}
\end{remark}

The proof of purely singular continuous spectrum has two
ingredients. The first is a proof of absence of absolutely
continuous spectrum. This follows by results of Kotani \cite{Kot}
and is in fact valid for every $\omega\in \Omega$ by results of
Last and Simon \cite{LS}. Alternatively, this follows by Theorem
\ref{zeromeasure} (whose proof, however, uses Kotani theory
\cite{Kot}). The second ingredient is a proof of absence of
eigenvalues. This is based on the so-called Gordon argument going
back to \cite{Gor}. Various variants of this argument have been
used in the study of \eqref{family} (see \cite{Dam} for a recent
overview). We use it in the following form \cite{Dam, DP,Kam}.

\begin{lemma} Let $(\Omega,T)$ be a uniquely ergodic subshift over $A$. Let $(n_k)$ be
a sequence in $\N$ with $n_k \to  \infty$, $k \to  \infty$. Set
$$
\Omega(k)\equiv  \{ \omega \in \Omega : \omega (-n_k +l) = \omega(l) = \omega(n_k + l),
\; 0 \leq l \leq n_k-1\}.
$$
If $\limsup_{k\to\infty} \mu(\Omega(k))>0$, then $\mu$-almost
surely, $H_\omega$ has no eigenvalues.
\end{lemma}

The lemma reduces the proof of absence of eigenvalues to
establishing the occurrence of sufficiently many cubes. For
primitive substitutions, occurrence of many cubes follows from
occurrence of one word of the form $uuue$, where $e$ is the first
letter of $u$. This was shown by one of the authors in \cite{Dam5}
(see \cite{Dam3} as well). It turns out that this line of
reasoning can be carried over to minimal substitutions. Namely, we
have the following result.

\begin{lemma}\label{cubes} Let $\Oomega$ be a minimal
substitution dynamical system. Let $u\in \CalW$ be given starting
with $e\in C$ such that $uuue$ belongs to $\CalW$. Set $n_k\equiv
|S^k (u)|$. Then,  $\limsup_{k\to \infty} \mu (\Omega(n_k)) >0$.
\end{lemma}

\begin{proof} As already mentioned, the proof is modelled after
\cite{Dam5}. As $uuue$ occurs in $\CalW$, so does $S^k (u u u e)$
for $k\in \N$. Of course, $S^k(u)$ begins with $S^k (e)$. Thus,
each occurrence of $S^k (u u u e) = S^k (u) S^k (u) S^k(u) S^k(e)$
gives rise to $|S^k(e)|$ occurrences of cubes and  we infer

\begin{equation}\label{eeins} \mu (\Omega(n_k)) \geq \mu (\Omega_{S^k (u u u e)})\times
|S^k(e)|,
\end{equation}
where we set $\Omega_v\equiv \{\omega\in \Omega : \omega(1)\ldots
\omega(|v|) = v\}$ for $v\in\CalW$. By Proposition~\ref{growth},
there exist $\lambda,\rho>0$ and $\theta>0$ with
\begin{equation}\label{ezwei}
|S^k(u u u e)| \leq  \rho \,\theta^k \;\: \mbox{and} \;\:  \lambda \,\theta^k
\leq |S^k (e) |,
\end{equation}
for every $k \in \N$. By Corollary \ref{ue}, $\Oomega$ is uniquely
ergodic and therefore
$$
\mu \left( \Omega_{S^k (u u u e)} \right) = \lim_{|x|\to
\infty}\frac{ \#_{S^k (u u u e)}
  (x)}{|x|}.
$$
Moreover, by Theorem \ref{main}, $\Oomega$ is linearly repetitive
with some constant $C_{{\rm LR}}$. Combining these estimates, we
infer

\begin{eqnarray*}
\mu (\Omega(n_k)) &\geq& \lim_{|x|\to \infty} \frac{ \#_{S^k ( u u u e)} (x)}{|x|} |S^k(e)|\\
&\geq & \frac{1}{C_{{\rm LR}}  |S^k (u u u e)|} |S^k (e)|\\
&\geq & \frac{\lambda}{C_{{\rm LR}}\, \rho}.
\end{eqnarray*}
This finishes the proof.
\end{proof}

\begin{proof}[Proof of Theorem~\ref{singularcont}.]
By the discussion following the theorem, it suffices to show
almost sure absence of point spectrum. This is an immediate
consequence of the preceding two lemmas.
\end{proof}

\begin{remark}
{\rm There is another approach to proving absence of eigenvalues
which is based on palindromes, rather than cubes. Concretely,
$\Omega$ is said to be palindromic if $\mathcal{W} (\Omega)$
contains arbitrarily long palindromes. Hof et al.\ prove in
\cite{HKS} that if $\Omega$ is minimal and palindromic, then for a
dense $G_\delta$-set of $\omega \in \Omega$, the operator
$H_\omega$ has empty point spectrum. This gives another method to
prove absence of eigenvalues for minimal substitution Hamiltonians
in cases where Lemma~\ref{cubes} does not apply, but where
sufficiently many palindromes occur. }
\end{remark}

\section{Fixed Points of Linearly Repetitive Substitutions} \label{Fixed}

In this section we discuss an application of the results above to number theory.

Recall that every $z\in (0,1)$ has a binary expansion
$$
z = \sum_{n=1}^\infty \frac{a_n}{2^n}
$$
with $a_n \in \{0,1\}$. For algebraic numbers, this binary
expansion is expected to be a ``random sequence.''  Of course,
there are various ways to give a precise meaning to ``random
sequence.'' One particular way is  that this binary expansion
should not be a fixed point of a substitution (see \cite{AZ} for
further discussion). The question whether such a binary expansion
can be a fixed point has thus attracted attention, and the most
general result so far has been obtained by Allouche and Zamboni
\cite{AZ}. Namely, they show that a binary expansion which is a
fixed point of either a primitive substitution or a substitution
of constant length (i.e., the images of letters all have equal
length) can only belong to a rational or transcendental number. We
can prove the following result, merely assuming minimality:

\begin{theorem}\label{number}
Suppose $S : \{ 0,1 \} \rightarrow \{ 0,1 \}^*$ satisfies {\rm
(BG)} {\rm (}i.e., $S$ induces a minimal, linearly repetitive
dynamical system{\rm )}. If $u \in \{0,1\}^\N$ is an aperiodic
fixed point of $S$ and $z \in (0,1)$ is given by
$$
z = \sum_{n=1}^\infty \frac{u_n}{2^n},
$$
then $z$ is transcendental.
\end{theorem}

\begin{proof}
For primitive $S$, the assertion was shown in \cite{AZ}. Let us
consider the case where $S$ is nonprimitive. Then the alphabet $B$
is not empty and we have either $S(0) = 0$ or $S(1) = 1$. Let us
discuss the case $S(1) = 1$, the other case can be treated in an
analogous way. From (BG) we can infer that $S(0)$ contains both
$0$ and $1$ and it begins and ends with $0$. By aperiodicity,
$S(0)$ cannot be equal to $0 1^k 0$ with $k\geq 1$.  That is,
either $S(0)$ has the form

\begin{equation}\label{firstcase}
S(0) = 0 1^k 0 w 0
\end{equation}
for some suitable word $w$, possibly empty, and suitable $k \ge 1$, or $S(0)$ has the form

\begin{equation}\label{secondcase}
S(0) = 00w0
\end{equation}
for some word $w$ containing $1$ as a factor.

We first consider the case where $S(0)$ is given by \eqref{firstcase}. Then
$$
S^2(0) = 0 1^k 0 w 0 1^k 0 1^k 0 w 0 S(w) 0 1^k 0 w 0,
$$
where we let $S(\varepsilon) = \varepsilon$ for definiteness. We
see that $u$ contains the word $0 1^k 0 1^k 0$ and hence, for some
prefix $p$,
$$
u = p 0 1^k 0 1^k 0 \ldots
$$
Now define
$$
U_n = S^n(p), \; V_n = S^n(0 1^k), \; V_n' = S^n(0).
$$
Observe that we have
\begin{equation}\label{fmeins}
|V_n| \rightarrow \infty \mbox{ as } n \rightarrow \infty
\end{equation}
by (BG),
\begin{equation}\label{fmzwei}
\limsup_{n \rightarrow \infty} \frac{|U_n|}{|V_n|} < \infty
\end{equation}
by Proposition~\ref{growth}, and
\begin{equation}\label{fmdrei}
\liminf_{n \rightarrow \infty} \frac{|V_n'|}{|V_n|} > 0,
\end{equation}
again by Proposition~\ref{growth}. We can now conclude the proof
in this case by applying \cite[Proposition~1]{FM} since
\eqref{fmeins}--\eqref{fmdrei} provide exactly the necessary input
for an application of this proposition.

Let us now consider the case where $S(0)$ is given by
\eqref{secondcase}. Then $S^2(0)$, and hence $u$, contains the
factor $0^3$. Therefore, $u = p000\ldots$, so we can set $U_n =
S^n(p)$, $V_n = V_n' = S^n(0)$ and then conclude as above.
\end{proof}

By the same argument one can prove the following extension; see
the note added in proof of \cite{AZ} for the necessary additional
input (namely, a result of Mahler \cite{M}).

\begin{theorem}\label{number2}
If $z \in (0,1)$ has a base $b$ expansion {\rm (}$b \in \N$ and $>
1${\rm )} which is given by an aperiodic fixed point of a
substitution on a two-letter alphabet which satisfies {\rm (BG)},
then $z$ is transcendental.
\end{theorem}

\section{Unique Decomposition Property}\label{Unique}
In this section we study questions concerning unique decomposition
for nonprimitive minimal aperiodic substitutions. For primitive
substitutions, such a unique decomposition property has been shown
by Moss\'{e} \cite{Mos} (for a study of the higher-dimensional
case, we refer to \cite{Sol}). We are not able to treat the
general case but rather restrict our attention to a two-letter
alphabet. This case has attracted particular attention recently in
the work of de Oliveira and Lima \cite{OL}.  Thus, we can assume
(and will assume throughout this section) that $A=\{a,b\}$ and
$$
|S(a)|>1 \quad \text{and} \quad S(b) = b.
$$
For $w\in \CalW$, an equation
$$
w = z_0 z_1\ldots z_n z_{n+1}
$$
is called a $1$-partition if $z_1 , \ldots , z_n \in \{S(a), b\}$,
$z_0$ is a suffix of an element in $\{S(a),b\}$ and $z_{n+1}$ is a
prefix of an element in $\{S(a),b\}$.  Similarly,  a $1$-partition
of $\omega\in \OmegaS$ consists of a sequence $(E_n)\subset \Z$
with
$$
\cdots < E_{-2} < E_{-1} < E_0 < E_1 < E_2 < \cdots, \;\:
\mbox{and}\:\; \lim_{n\to  \infty} E_n= \infty,   \;\:\lim_{n\to
-\infty} E_n= -\infty,
$$
such that $\omega(E_n) \ldots \omega (E_{n+1} -1) \in \{S(a),b\}$.
We will prove the following theorem.

\begin{theorem}\label{decomposition}
Let $\Oomega$ be a minimal, aperiodic, nonprimitive substitution
dynamical system over $\{a,b\}$.  Then, every $\omega\in \OmegaS$
admits a unique $1$-partition. More precisely, there exist $L\in
\N$ and words $w_1,\ldots, w_k$ such that the $1$-partition
$\{E_j\}$ of $\omega$ consists of exactly those $E$ with $\omega
(E-L)\ldots \omega (E +L) \in \{w_j : j=1,\ldots ,k\}$.
\end{theorem}

To prove this result, we need some preparation. The proof of
Theorem~\ref{decomposition} appears at the end of this section.

\begin{lemma}\label{basicR}
Let $\Oomega$ be a minimal, aperiodic, nonprimitive substitution
dynamical system over $\{a,b\}$. Then $S(a)$ is neither a prefix
of $b S(a)$ nor a suffix of $S(a) b$. In particular, two
$1$-partitions of a finite $v$ which both start with $S(a)$ {\rm
(}end with $S(a)${\rm )} must agree up to a suffix {\rm
(}prefix{\rm )} of length at most $ |S(a)b|$.
\end{lemma}

\begin{proof}
The second statement follows immediately from the first. So assume
the first statement is wrong. Then, $S(a) = b^l$ with a suitable
$l\in \N$ and periodicity of $\Oomega$ follows.
\end{proof}

\begin{lemma}\label{l62}
Let $\Oomega$ be a minimal, aperiodic, nonprimitive substitution
dynamical system over $\{a,b\}$. Set
$$
L_0 \equiv \max\{ |v^n| : v \;\:\mbox{subword of $S(a)$ and $v^n
\in \CalW$}\}.
$$
If $v$ with $|v|>2 |S(a)| + L_0$ admits a $1$-partition beginning
with $S(a) S(a)$, then every $1$-partition of $v$ starts with
$S(a)S(a)$.
\end{lemma}

\begin{proof}
Assume the contrary. By Lemma~\ref{basicR}, there exists then a
$1$-partition $v = z_0 z_1 z_2 \ldots z_n z_{n+1}$ of $v$ with $0
< |z_0| < |S(a)|$ and $z_1 = S(a)$. This gives $S(a) = v^r$ with a
primitive $v$ and $r\in\N$ suitable. By definition of $L_0$, both
this $1$-partition of $v$ and the $1$-partition beginning with
$S(a) S(a)$ contain blocks of the form $b$. Consider the leftmost
of these blocks. Then, we obtain that $v$ is a suffix of $v b$ and
thus, $v=b^l$. This in turn yields the contradiction $S(a) = b^{r
|b| }$.
\end{proof}

It will be convenient to treat the two cases, where $S(a)$ does or
does not contain the word $aa$, separately. We first consider the
case where $aa$ is a subword of $S(a)$ and prove uniqueness of
decompositions under this assumption.

\begin{prop}\label{aaprop}
Let $\Oomega$ be a minimal, aperiodic, nonprimitive substitution
dynamical system over $\{a,b\}$. Suppose that $aa$ occurs in
$S(a)$. Then, there exists an $L\in \N$ such that all
$1$-partitions of $v\in \CalW$ induce the same $1$-partition on
$v(L) \ldots v (|v| -L)$.
\end{prop}

\begin{proof}
By Lemma~\ref{l62}, all occurrences of the word $S(a)S(a)$ in
$1$-partitions of $v$ which begin before the $L_0$-th position in
$v$ are uniquely determined. In particular, they occur at the same
places in all $1$-partitions. By Lemma~\ref{basicR}, the
$1$-partitions must then agree to the right and to the left of
such an occurrence up to a boundary term of length not exceeding
$|S(a)b|$.  Thus, it suffices to show existence of a $1$-partition
of $v$ containing the blocks $S(a) S(a)$ if $v$ is long enough.
This, however, is clear by the assumption that $aa$ occurs in
$S(a)$ and minimality of $\Oomega$, as every $v$ is contained in a
word of the form $S^n (a)= S (S^{n-1} (a))$.
\end{proof}

Next we turn to the case where $S(a)$ does not contain $aa$ as a
factor.

\begin{lemma}\label{rgleichr}
Let $\Oomega$ be a minimal, aperiodic, nonprimitive substitution
dynamical system over $\{a,b\}$. Suppose that $aa$ does not occur
in $S(a)$. Let $w\in \CalW$ be given with
$$
w = S(a) b^{r_1} S(a) b^{r_2} S(a) \cdots S(a) b^{r_n} S(a) = x
S(a) b^{s_1} S(a) b^{s_2} S(a) \cdots S(a) b^{s_u} y
$$
with suitable $r_1,\ldots, r_n\in \N$, $s_1,\ldots,s_u\in \N$ and
$x,y\in \CalW$ with $|x| < |S(a)|$ and $|y|\leq |S(a)|$. Then,
$u=n$ and $r_1 = s_1 = r_2 = s_2 = \cdots = s_n = r_n$.
\end{lemma}

\begin{proof}
By minimality, $S(a)$ begins and ends with $a$. As
$a a$ does not occur in $S(a)$, we have
$$
S(a) = a b^{k_1} a \ldots a b^{k_l} a
$$
with suitable $k_1, \ldots, k_l \in \N$. If $|x|=0$, the statement
now follows easily. Otherwise, we have $x = a b^{k_1}a \ldots a
b^{k_j}$ with $j< l$ suitable. Consider the blocks of consecutive
$b's$ appearing in $w$. Such a block will be referred to as
$b$-block. By $w= S(a) b^{r_1}\ldots$, we infer that $b^{r_1}$ is
the $(l+1)$-st $b$-block appearing in $w$. By $w= x S(a)\ldots$,
we then infer that $b^{r_1}$ occurs in $S(a)$ as the $(l+ 1 -
j)$-th $b$-block. Similarly, considering the occurrence of
$b^{r_2}$ in $w = S(a) b^{r_1} S(a) b^{r_2}\ldots$, we infer that
$b^{r_2}$ is the $(2 l +2)$-th $b$-block appearing in $w$. On the
other hand, by $w= x S(a) b^{s_1} S(a) \ldots$, we see that
$b^{r_2}$ is the  $j + l + 1 + t$-th $b$-block in $w$, where $t-1$
is the relative number of $b$-blocks occurring in $S(a)$
before $b^{r_2}$. This yields
$$ 2 l + 2 = j + l +1 + t,$$
and we infer $t = l +1 -j$. Thus, $b^{r_1}$ and $b^{r_2}$ occur in
the corresponding $S(a)$ blocks at the same relative positions.
This yields immediately $r_1=r_2$.

Denote the relative position of $b^{r_1} = b^{r_2}$ in $S(a)$ by
$p$. Thus, $p= |S(a)| - |x| + 1$ by $w= S(a) b^{r_1} \ldots = x
S(a) \ldots$. Now, consider the absolute position $h$ of $b^{r_2}$
in $w$. Then,
$$ h = 2 |S(a)| + r_1 + 1  = |x| + |S(a)| + s_1 + p.$$
Putting this together, we infer $r_1 = s_1$. Now, the assertion
follows easily by repeating this reasoning.
\end{proof}

\begin{prop}\label{aprop}
Let $\Oomega$ be a minimal, aperiodic, nonprimitive substitution
dynamical system over $\{a,b\}$. Suppose that $aa$ does not occur
in $S(a)$. Then, there exists $L\in \N$ such that all
$1$-partitions of $v\in \CalW$ induce the same $1$-partition on
$v(L) \ldots v (|v| -L)$.
\end{prop}
\begin{proof} As $S$ is minimal and aperiodic, there exists an
$N\in \N$ such that $v^n\in \CalW$ for $v\in \CalW$ with $|v|\leq
2 |S(a)|$ implies $n\leq N$. Set $L = (N +2) 2 |S(a)|$. Assume
that there exists a $w\in \CalW$ admitting two $1$-partitions
inducing two different $1$-partitions of $w(L)\ldots w(|w|-L)$. By
Lemma~\ref{rgleichr}, we infer existence of an $(N+1)$-power of
$v= S(a) b^r$ in $\CalW$. This contradicts the choice of $L$.
\end{proof}

\begin{proof}[Proof of Theorem \ref{decomposition}.]
Let $\omega\in \Omega$ be given. Uniqueness of the $1$-partition
is clear from Propositions~\ref{aaprop} and \ref{aprop}. Existence
of a $1$-partition follows by standard compactness-type arguments
but can also be shown as follows: For $n\in \N$, set $v_n \equiv
\omega(-n)\ldots \omega (n)$. Then each $v_n$ admits a
$1$-partition as it is contained in $S^l (a)$ with a suitable $l$.
These $1$-partitions are compatible due to the previous
proposition. Thus, they easily induce a $1$-partition of $\omega$.
These considerations and the previous proposition also imply the
last statement of the theorem.
\end{proof}

{\it Acknowledgements.} The authors would like to thank 
Hisatoshi Yuasa  for pointing out a mistake in an earlier version of the manuscript and providing the example discussed in Remark 1.b.

\end{document}